\documentclass[12pt, twoside, a4paper, reqno]{amsart}
\usepackage{amsmath, amssymb, latexsym, amsthm, amsfonts, mathrsfs, amscd}
\DeclareMathAlphabet{\mathpzc}{OT1}{pzc}{m}{it}
\usepackage{enumerate}[1.)]
\renewcommand{\eqref}[1]{(\ref{#1})}   
\usepackage{latexsym}
\usepackage{euscript}
\usepackage[margin=2.95 cm,bmargin=3.5 cm,tmargin=3.5 cm]{geometry}
\usepackage{epsfig}
\usepackage{hyperref}
\usepackage{url}
\usepackage{tikz}
\usepackage{graphics}
\usepackage{array}
\usepackage{enumerate}

\usepackage[utf8]{inputenc}
\usepackage{xcolor}
\newtheorem{theorem}{Theorem}[section]
\newtheorem{lemma}[theorem]{Lemma}
\newtheorem{proposition}[theorem]{Proposition}
\newtheorem{corollary}[theorem]{Corollary}

\newtheorem{remark}{Remark}[section]

\renewcommand{\l}{\left}
\renewcommand{\r}{\right}

\newcommand{\moda}{(\textnormal{mod }a)}

\newcommand{\diff}{\textnormal{d}}

\newcommand{\R}{{\mathbb R}}
\newcommand{\Z}{{\mathbb Z}}
\newcommand{\C}{{\mathbb C}}

\newcommand{\abcdsumX}{\mathop{\sum_{p_1}\sum_b\sum_{p_2}\sum_d}_{\substack{p_1d-p_2d=r\\ -X \leq p_1, b, p_2, d \leq X \\ (m_1m_2p_1p_1,q) = 1 }}}

\newcommand{\be}{\begin{equation}}
\newcommand{\ee}{\end{equation}}
\newcommand{\ba}{\begin{equation}\begin{aligned}}
\newcommand{\ea}{\end{aligned}\end{equation}}
\numberwithin{equation}{section}
 2

\usepackage{setspace}
\usepackage{appendix}

\newcommand{\ve}{\varepsilon}
\numberwithin{equation}{section}

\begin{document}
\title[Prime entries]{Prime Solutions to a Binary Additive Equation and Mixed Moments of Character Sums}  
\author{Rachita Guria}
\address{Max Planck Institute for Mathematics, Vivatsgasse 7, 53111, Bonn, Germany}
\email{guriarachita2011@gmail.com}

\maketitle	
	
\begin{abstract}
We obtain an asymptotic formula with a power-saving error term for counting the integer points $(a,b,c,d)$ in an expanding box that satisfy the determinant equation $x_1x_2-x_3x_4 =r$ for $r \neq 0 $ with two of entries to be prime. Finally, these estimates are applied to evaluate mixed fourth moments of Dirichlet character sums over integers and primes, yielding non-trivial bounds. The method involves the Poisson summation formula and the estimation for the average of the sums of the Kloosterman fractions over primes.
\end{abstract}
\section{Introduction}

The Waring–Goldbach problem lies at the intersection of two classical questions in number theory: Waring’s problem and Goldbach’s conjecture. It asks whether large integers can be expressed as sums of a bounded number of prime powers. A notable case is the classical Lagrange-type equation:
\ba \label{lagrange}
x_1^2+x_2^2+x_3^2+x_4^2 = N
\ea
with multiplicative restrictions, such as requiring the variables to be prime.. It is plausible to expect that sufficiently large integers, under certain necessary congruence conditions, can be written as the sum of four squares of primes. Such a result of this strength remains out of reach and is still a conjecture. However, Hua in \cite{hua1947additive} proved that all large integers congruent to $5$
modulo $24$ can be represented as the sum of five squares of primes. This was achieved by the Vinogradov’s method, originally used for the ternary Goldbach problem which is a ternary additive problem like the representatin of an integer as a sum of five prime squares. For the solutions to ternary additive problems like $\alpha+\beta+\gamma = N$, we expect a non-trivial cancellation through the circle method since at least two of the sequences are sufficiently dense in $[1,N]$. But binary problems lack this advantage.\\

In this article, we are interested in the binary additive problem. A relevant binary additive problem is the solutions to the Lagrange equation with almost prime variables.. In 1994, Brüdern and Fouvry in \cite{brudern1994lagrange} established that for every sufficiently large integer $N$, congruent to $4$ modulo $24$, \eqref{lagrange} admits solutions over integers, each of which is almost prime of order $34$ by combining a sieve method with the Hardy-Littlewood-Kloosterman approach. Concerning the Lagrange equation with four almost prime variables, the value $34$ due to Brüdern and Fouvry was later sharpened by Heath-Brown and Tolev in \cite{heath2003lagranges} to $25$, by Tolev in \cite{tolev2002lagrange} to $21$, by Cai in \cite{cai2010lagrange} to $13$ and by Tsang and Zhao in \cite{tsang2017lagrange} to $4$ .\\

In 2003, Heath-Brown and Tolev 
in \cite{heath2003lagranges} proved that equation \eqref{lagrange} has solutions with one prime entry and each of the other three $x_2,x_3,x_4$ having at most $101$ prime divisors. Tolev in \cite{tolev2002lagrange} reduced the number of prime divisors from $101$ to $80$, and then Cai in \cite{cai2010lagrange} showed that $42$ is acceptable. The latest improvement has been made in 2017 by Tsang and Zhao in \cite{tsang2017lagrange}. \\

Many other variants of the four squares theorem (see \eqref{lagrange}) have been studied as an approximation to the conjecture above. Podsypanin in \cite{podsypanin1975representation} obtained an asymptotic formula as a function of $N$ for the number of solutions of \eqref{lagrange} in square-free numbers $x_i$. \\

Greaves in \cite{greaves1976representation} applied a sieve method to solve \eqref{lagrange} with $x_1$, $x_2$ as integers and $x_3$, $x_4$ as primes, obtaining only a lower bound for the number of such representations. Shields in \cite{shields1979some}, Plaksin in \cite{plaksin1982asymptotic} and Kovalchik in \cite{koval1984some}, by different methods, obtained an asymptotic formula for such representations where the authors use the multiplicative properties of the arithmetic function $r(n)$ defined as the number of ways $n$ can be represented as the sum of two squares along with the methods of Hooley and Linnik's dispersion method. The best known result today is due to Plaksin in \cite{plaksin1985asymptotic} who obtained an error term of size $O\l(N(\log N)^{-2.042} \r)$ by proving the following theorem under some assumptions.
\begin{theorem}(\textbf{Plaksin})
Let $\kappa(N)$ be the number of solutions of the equation
\ba 
f (p,q) + \psi(x,y) = N,
\ea 
and $\tau(m,N)$ be  the number of solutions of the congruence
\ba \label{con}
f(u,v) + \psi(x,y) \equiv N \textnormal{ mod }m,\quad \quad (uv,m) = 1,
\ea
where $f$ and $\psi$ are positive definite quadratic forms with integral coefficients of discriminants $-\Delta$ and $-\tau$ respectively, $p$ and $q$ are prime variables, and $x,y,u$ and $v$ are integers. We let $p^{\beta}||(\Delta, 4N)$ and $G = \prod_{p | 3\Delta N, \tau} p ^{\beta +1}$. Let the congruence \eqref{con} be solvable for $m = G$. Then we have 
\ba 
\kappa(N) = \frac{8 \pi}{\sqrt{\tau}} s_{f} \mathfrak{S}(N)  \frac{N}{\log ^2 N}\l(1+ \frac{B}{\log ^{0.042}N}\r) ,
\ea
where the singular series has the form
\ba
\mathfrak{S}(N) = \prod_p \frac{\tau(p^{\beta +1},N)}{(\phi(p^{\beta+1}))^2p^{\beta+1}},\quad \textnormal{ and } \quad \frac{1}{\sqrt{G} \log\log N} \ll \mathfrak{S}(N) \ll \log \log N ,
\ea
and $s_{f}$ is the area of the region $f(x,y) \leq 1$ such that $x,y \geq 0$ and $\phi$ is the Euler totient function.
\end{theorem}

Our goal is to study a similar type of problem, specifically on the determinant equation: 
\ba \label{determinant}
x_1x_2-x_3x_4 =r,
\ea
with multiplicative constraints imposed on the variables and a non-zero $r$, instead of the Lagrange equation \eqref{lagrange}. It is important to note that solutions to both Lagrange’s equation and the determinant equation cannot be considered ternary additive problems, as their terms do not combine in a way that meets the relevant criteria (see \cite[Eq. (0.1.6) and (0.1.9)]{linnik1963dispersion}). However, they fit as examples of the binary additive problems as discussed earlier.\\

We focus on counting integer points in an expanding box $[-X,X]^4$ that satisfy \eqref{determinant}, where one of $x_3$ or $x_4$ is prime, and an arbitrary weight $\alpha(x_1)$ is attached to one variable, say $x_1$. Precisely, Our aim is to derive an exact asymptotic formula with a strong error bound, while allowing $r$ to grow as large as possible relative to $X$. This seems to be the first instance of an asymptotic formula with a power saving error term for counting solutions to a quaternary equation with two prime variables.\\

Exploiting the special structure of the determinant form, we show the following.


\begin{theorem}\label{mainthm}
For an arbitrary sequence of complex numbers $\alpha(n)=O(n^{\ve})$, and any non-zero integer $r$, let us define 
\ba 
\mathpzc{S}_r(X) = \mathop{\sum\sum\sum\sum}_{\substack{  |a|,|b|,|p|,|d| \leq X \\ ad- pb =r}} \alpha(a),
\ea
where $a,b,d$ vary over the integers and $p$ varies over primes. Then, for fixed $r$, we have,
\ba 
\mathpzc{S}_r(X) &= 4\mathop{\sum_{1 \leq a \leq X}\sum_{1 \leq p \leq X}} \frac {\alpha(a)} {a}  \int\limits _{1}^X \l[w\l(\frac{|r| +px}{a} \r) \diff x + w\l(\frac{-|r|+px}{a} \r) + w\l(\frac{|r| -px}{a} \r)\r] \diff x \\
&+ O\l(X^{1+\frac{3}{4}+\ve}+ r^{\frac{1}{5}}X^{1+\frac{11}{20}+\ve} \r),
\ea
as $X \rightarrow \infty$, and the weight function $w$ satisfies the conditions \eqref{w0}-\eqref{w-prime}.
\end{theorem} 
From the properties (see \eqref{w0}-\eqref{w-prime}) of the weight function $w$ it is easy to see that for $|r| \gg X^2$ the main term vanishes, and in the complementary case the main term is of size $X^{2+\ve}$. \\
This theorem gives an asymptotic uniform formula for $\mathpzc{S}_r(X)$ when $0 < |r| \ll X^{2-\ve}$. \\
If we take $\alpha(n)$ to be the characteristic function of primes, then we have the following theorem which gives us an asymptotic formula regarding counting lattice points with two prime entries on the determinant surfaces. We consider a more general sum.
\begin{theorem} Let q be an ineteger. \label{theorem}
Define
\be
\mathpzc{S}(2,\Z):=  \l\{\begin{pmatrix} p_1 & b\\ p_2 &
  d\end{pmatrix}\in M(2,\mathbb{Z}): p_1,p_2 \textnormal{ are primes and } (m_1m_2p_1p_2,q) = 1  \r\} . \nonumber
  \ee 
Let
\[
\tilde{\mathpzc{S}}_r(X) = \sum_{\substack{\gamma \in  \mathpzc{S}(2, \Z)\\ det(\gamma)=r\\
\\ ||\gamma||_{\infty} \leq X}} 1=\abcdsumX 1,
\]
where $r$ is a fixed integer. Then, for large $X$ and $r \neq 0$, we have
\ba
\tilde{\mathpzc{S}}_r(X)  &=  \frac{ 4\phi(q)}{q}\mathop{\sum_{1 \leq p_1 \leq X}\sum_{1 \leq p_2 \leq X}} \frac {1} {p_1} \int\limits _{1}^X \l[w\l(\frac{|r| +p_2x}{p_1} \r) \diff x + w\l(\frac{-|r|+p_2x}{p_1} \r) + w\l(\frac{|r| -p_2x}{p_1} \r)\r] \diff x \\
 &+ O\l(X^{1+\frac{3}{4}+\ve}+ r^{\frac{1}{5}}X^{1+\frac{11}{20}+\ve} \r).
\ea
and for $r = 0$, we have 
\[ \tilde{\mathpzc{S}}_0(X) = \frac{4\pi^2}{3}
\left(\prod_{p\mid q}\left(1-\frac{1}{p^2}\right)\right)
\frac{X^2}{\log^2 X}  + \frac{\phi(q)}{q} \pi(X) X 
+ O \left(\frac{X^2}{\log^3 X}\right).
\]
\end{theorem} 
\begin{remark}
The important thing to notice here is that even if we have restricted two of entries of the determinant equation (see \eqref{determinant}) to be prime, we still have been able to achieve a power saving in the error term, in contrast, for the Lagrange equation (see \eqref{lagrange}) with two prime entries only the logarithmic saving of size $(\log N)^{2.042}$ in the error term has been obtained so far.
\end{remark}

\begin{remark}
Following the same kind of initial reduction as in \cite[Eq. (3.2), \S 3]{GG} and the main term calculation from the same paper (see \cite[\S 4]{GG}) we can explicitly evaluate the main term in our case. For $r \neq 0 $, by Prime Number Theorem we show that
\ba 
\tilde{\mathpzc{S}}_r(X) &=  \frac{ 16\phi(q)}{q} \int_2^X \frac{\text{li}(t)}{t \log t} \, dt + O\left( X \exp(-c \sqrt{\log X}) \right) \\
&+ O \l(X^{\ve}\l( |r| +  X^{1+\frac{3}{4}}+ r^{\frac{1}{5}}X^{1+\frac{11}{20}}\r) \r),
\ea
where 
\ba
\text{li}(x) = \int_2^x \frac{dt}{\log t}.
\ea 
Note that, in this derivation, the contribution of $S^+_r(X)$ (see \eqref{decom1}) i.e., the third main term in the previous expression is absorbed into the error term by $O\l(r^{1+\ve}\r)$. For small $r$, this main term is more explicit and does not depend on the weight function.
\end{remark}

In \cite[Theorem 1.2]{GG}, Ganguly and Guria have obtained an asymptotic formula for one prime entry with an error term of size $O\l(X^{5/3+\ve}\r)$. In \cite{BC} and before that in \cite{DFI2}, sums of the form
\ba\label{bc}
\mathop{\sum_{n_1}\sum_{n_2}\sum_{n_3}\sum_{n_4}}_{\substack{n_i\in \mathcal{N}_i 1\leq i\leq 4\\n_1 n_4-n_2n_3=r }}\alpha(n_1)\beta(n_2)f(n_3)g(n_4)
\ea
were studied and asymptotic formulae for such sums were obtained. Here, $\mathcal{N}_i=[N_i/2, N_i]$, $\alpha$ and $\beta$ are arbitrary sequences, and $f, g$ are smooth functions satisfying suitable decay conditions on their derivatives. One may compare our theorems with \cite[Corollary 1]{BC} or \cite[Theorem 1]{DFI-det}, where, two sequences are arbitrary and the other two are smooth weights of the form $f(n_3)$ and $g(n_4)$ (as in \eqref{bc}). Note that if we apply these results in our situations, we obtain an error term that is rather poor. For example, \cite[Corollary 1]{BC} gives an error term of the order $O\l(X^{1+49/50+\ve}\r)$. \\

The following two corollaries are direct consequences of Theorem \ref{mainthm}.
\begin{corollary}
We have,
\ba 
\mathop{\sum\sum\sum}_{\substack{1 \leq a,b,c \leq X }} \Lambda\l(\frac{ab\pm r}{c} \r) \alpha(a) &= \sum_{ p } \Lambda(p)  \sum_{ 1 \leq a \leq X} \frac{\alpha(a)}{a} \int\limits_1^X w\l(\frac{xp\mp r}{a} \r) \\
&+ O\l(X^{1+\frac{3}{4}+\ve}+ r^{\frac{1}{5}}X^{1+\frac{11}{20}+\ve} \r),
\ea
uniformly for all non-zero integer $r$ and $a,b,c$ vary over the integers, and the weight function $w$ satisfies the conditions \eqref{w0}-\eqref{w-prime}.
\end{corollary}
\begin{corollary} For any non-zero integer $r$, we have
\ba 
\mathop{\sum\sum\sum}_{\substack{1 \leq a,b,p \leq X }} \alpha \l(\frac{ap\pm r}{b} \r) &= \sum_{1 \leq p \leq X}\sum_{n \in \Z} \frac{\alpha(n)}{n} \int\limits_1^X w\l(\frac{xp\pm r}{n} \r) \\
&+ O\l(X^{1+\frac{3}{4}+\ve}+ r^{\frac{1}{5}}X^{1+\frac{11}{20}+\ve} \r),
\ea
where $a,b$ vary over the integers and $p$ varies over primes, and the function $w$ satisfies the conditions \eqref{w0}-\eqref{w-prime}.
\end{corollary}
\subsection*{Mixed moments of Dirichlet character sums} Our second question concerns the mixed fourth moment of character sums, where one of the sums is over integers and the other is over prime numbers. While classical mean value theorems typically focus on either on integers or on primes, we investigate the arithmetic behavior governing the interaction between these two distinct structures. Specifically, we analyze the following sum: 
\ba 
\mathcal{M}(X,q) = \sum_{\chi (\bmod q)} \left| \sum_{m \leq X}  \chi(m)  \right|^2 \left| \sum_{p \leq X}  \chi(p)  \right|^2 .
\ea 
Our initial motivation for investigating such sums stems from their connection to the Generalized Riemann Hypothesis (GRH). A central challenge in number theory is understanding the cancellation within character sums over primes, $\sum_{p\leq X} \chi(p)$. Under GRH, these sums exhibit optimal square-root cancellation. By embedding this prime sum into a family-wide average, we apply our matrix estimates to unconditionally evaluate $\mathcal{M}(X,q)$. Such weighted sums do occur naturally in practice, e.g., while studying the subconvexity problem for Dirichlet $L$-functions. The amplification method, a framework that requires evaluating mean values where a smooth sum over integers is inherently coupled with a sparse sum over primes, is a well known approach to subconvexity problem. Our mixed moment represents the structural core of such method. We now state several more results that demonstrate the applicability of our previous results.
\begin{theorem}\label{mixedmoment} We have
\ba 
\sum_{\chi (\bmod q)} \left| \sum_{m \leq X} \sum_{p \leq X} \chi(mp)  \right|^2 &=
\sum_{\chi (\bmod q)} \left| \sum_{m \leq X}  \chi(m)  \right|^2 \left| \sum_{p \leq X}  \chi(p)  \right|^2  \\
&= \l( \frac{\phi(q)}{q}\r)^2 X^2\pi(X)^2
+ O \left(X^{3+19/20+\ve} +qX^{2+\ve} \right),
\ea 
where $p$ runs over primes.
\end{theorem}
By isolating the contribution of the principal character $\chi_0$ in the sum on the left-hand side, the main terms cancel out, yielding the following corollary.
\begin{corollary}
We have
\ba  
\sum_{\chi \neq \chi_0(\bmod q)} \left| \sum_{m \leq X}  \chi(m)  \right|^2 \left| \sum_{p \leq X}  \chi(p)  \right|^2  \ll X^{3+19/20+\ve} +qX^{2+\ve} .
\ea
\end{corollary}
\begin{remark}
We obtain a non-trivial bound with X as large as $q^{20/39}$ which is greater than $\frac 1 2$. In \cite{DFI2}, the exponent was $q^{49/95}$.
\end{remark}
These types of sums were previously studied by Ayyad, Cochrane, Zheng \cite{ACZ}, by Cochrane, Sih \cite{cochone} and Kerr \cite{kerr}, in the case when the summations are taken over the integers.
\subsection{Method of the proofs}
We briefly outline the method of the proof of Theorem \ref{mainthm} since the remaining theorems are deduced from this result. The main idea is to exploit the special structure of the determinant equation \eqref{determinant} that enables us to transform our sum into the average of the Kloosterman fractions over primes $p$ after using the Poisson summation formula.\\

At first, we isolate terms for which $p|a$ and bound them trivially. Now, we are left with the terms which satisfy the condition $(a,p)=1$. The basic idea is to use Fourier analysis. To do so we first approximate this sum by a weighted sum by introducing a smooth function $w$ in the those variables without any restriction, namely $b,d$, which runs over integers, that approximates the indicator function $\mathbf{1}_{[1,X]}$ to within an acceptable error term (see \S \ref{weightfunction}). Now, we replace $d$ by $(r+pb)/a$, using the determinant equation and the congruence condition $b \equiv -r\bar{p} \moda$ arises as we already have the condition $(a,p)$ =1. At this stage, a trivial estimation gives $S_w(X,r) \ll X^{2+\ve}$. Next, we use the Poisson summation
formula to evaluate the sum over $b$ since it is the only variable that runs over all integers. We separate the zero frequency and the non-zero frequencies. The zero frequency yields the main term, whereas in regard to the non-zero frequencies, the key observation is that the sum over primes $p$ along with sum over modulus $a$ is precisely the average of the Kloosterman fractions over primes. Due to cancellations in the sum of Kloosterman fractions over primes, we expect the sum of the nonzero frequencies to be small and this is the crucial step in this proof. Such sums came up in an earlier investigation by Fouvry and Shparlsinki \cite{FSternary} and applying their Theorem 3 one can get a power saving error term. However, this result has been improved by Irvine (see Lemma \ref{exponential}) and we apply this instead. It requires some work because of the condition $Q^{4/3} \geq X \geq Q^{1/2}$. We also need a good bound for our oscillatory integral (coming from $\l( \widehat{F_{da_1,t}}\left(\frac{n}{da_1}\right)  \r)'$) as the intervals over which the variables run are long and the direct application of integration by parts is not enough.
 \subsection*{Acknowledgements}
The author thanks Valentin Blomer, Satadal Ganguly, Pieter Moree, Alina Ostafe, Prahlad Sharma, and Igor Shparlinski for helpful comments. She is grateful to Max Planck Institute for Mathematics in Bonn where this work was carried out for its hospitality, excellent working atmosphere and financial support.

\section{Preliminaries}
In this section we recall some standard formulae and lemma from the literature which will be required in the proofs of our results.

\begin{lemma}\label{10}\textbf{(Application of Poisson summation)}\\
Suppose that both $f$, $\hat{f}$ are in $L^1(\R)$ and have bounded variation. Then 
\be \sum_{\substack{n \equiv \alpha \moda \\ n \in \Z}}f\left(n\right)= \frac{1}{a}\sum_{n \in \Z}e \left(\frac{\alpha n}{a}\right)\hat{f}\left(\frac{n}{a}\right) \ee where both the series converge absolutely.
\end{lemma}
\begin{proof}
See, e.g., \cite[Theorem 4.4 and Exercise 4]{IK}.
\end{proof}
\begin{lemma} \label{exponential}
For any integer $a>0$ and every $\ve > 0$, we have
\be 
\sum_{q \sim Q} \l| \sum_{\substack{ 1\leq p \leq X \\ (p,q)=1}} e\left(\frac{a \Bar{p}}{q}\right)\r| \ll
\l(1+\frac{a}{XQ}\r)^{\frac{1}{2}}   \l(Q ^{\frac{1}{2}}X^{\frac{11}{8}}+ Q^{\frac{7}{6}} X^{\frac{2}{3}}\r)(aQ)^{\ve},
\ee 
for $Q^{\frac{4}{3}} \geq X \geq Q^{\frac{1}{2}}$ .
\end{lemma}
\begin{proof}
See \cite[Theorem 1.4]{irving2014average} and \cite[Theorem 5]{FSternary} .
\end{proof}

\section{The main steps of proof of Theorem \ref{mainthm}}
Let us define
\ba \label{decom1}
S_r^-(X) := \mathop{\sum\sum\sum\sum}_{\substack{1 \leq a,b,p,d \leq X \\ ad- pb =r}} \alpha(a), \textnormal{ and } S_r^+(X) := \mathop{\sum\sum\sum\sum}_{\substack{1 \leq a,b,p,d \leq X \\ ad + pb =r}} \alpha(a),
\ea
and we have
\ba \label{thm}
\mathpzc{S}_r(X) = 4S_r^- (X) + 4S_{-r}^- (X) + 4 S_{|r|}^+(X) + O\l( X^{1+\ve }\r)
\ea
Note that, we cannot use the symmetry because of the presence of the factor $\alpha(a)$ in our sum. It is enough to focus on the sum $S_r^- (X)$ for $r > 0$ since the analysis of the other two sums can be done in a much of a similar manner. Hence, we start with 
\ba 
S_r^-(X) &= \mathop{\sum_{1 \leq a \leq X}\sum_{1 \leq p \leq X}} \frac {\alpha(a)} {a} \int\limits _{1}^X w\l(\frac{ r +px}{a} \r) \diff x  \\
&+ O\l(X^{1+\frac{3}{4}+\ve}+ r^{\frac{1}{5}}X^{1+\frac{11}{20}+\ve} \r).
\ea
\subsection{Contribution of the terms when $p|a$}
We will now examine the contribution of the terms for which $p|a$ so that later we will only focus on the case when $(a,p)=1$. Hence, we have the sum 
\ba 
\mathop{\sum\sum\sum\sum}_{\substack{1 \leq a,b,p,d \leq X \\ ad- pb =r \\ p|a}} \alpha(a) 
&= \mathop{\sum\sum\sum}_{\substack{1 \leq b,p,d \leq X \\ kpd- pb =r \\ p|r}}\sum_{1 \leq k \leq \frac{X}{p}} \alpha(pk) \\
& = \sum_{\substack{1 \leq p \leq X \\ p|r}} \sum_{1 \leq k \leq \frac{X}{p}} \alpha(pk) \sum_{ \substack{1 \leq b \leq X \\ pk|(bp+r)}} 1\\
& \ll X^{\ve}\sum_{\substack{1 \leq p \leq X \\ p|r}} \sum_{1 \leq k \leq \frac{X}{p}} \frac{X}{k}\\
&\ll  X^{1+\ve} ,
\ea
since there are $O\l(X^{\ve}\r)$ values of $p$ that contribute in the sum over $p$ and $\alpha(n)= O(X^{\ve})$. Hence, 
\ba \label{Srx}
S_r^-(X) = S_r^{(0)}(X) + O\l( X^{1+\ve}\r),
\ea
where 
\ba
S_r^{(0)}(X) =  \mathop{\sum\sum\sum\sum}_{\substack{1 \leq a,b,p,d \leq X \\ ad- pb =r \\ (a,p) =1}}\alpha(a)
\ea
We next introduce a nice suitable weight function $w$ to the variables $b$ and $d$ in our sum $S_r^{(0)}(X)$ and we use the same weight function as in \cite[\S 3]{GG}. Let us recall the weight function $w$. \\

\subsection{Introducing a weight function}\label{weightfunction} Let $w:\R\to \R$ be a smooth function satisfying the following conditions:
\begin{align}
&(1)  \ w(x)=0\textnormal{ if } x\not\in (1-H, X+H),\label{w0}\\
& (2) \  w(x)=1\textnormal{ if } x\in [1, X], \label{w1}\\
& (3) \ 0<w(x)<1 \textnormal{ if } x\in (1-H, 1), \textnormal{ or if } x\in (X, X+H),\label{w-support}\\
& (4)\ w^{(j)}(x) \ll_j H^{-j}\textnormal{ for every } j\geq 1 \label{w-prime},
\end{align}
where $X>0$ is a large number and $H$ is a parameter to be fixed later subject to the condition 
\be\label{H-size}
\sqrt{X}\leq H\leq X.
\ee
\begin{remark}
Finally, our choice will be  $H=\max\l\{X^{\frac{3}{4}}, r^{\frac{1}{5}}X^{\frac{11}{20}}\r\}$.
\end{remark}
We define next
\ba \label{defSrwX} 
S_w(X,r) &= \mathop{\sum_{ 1 \leq a \leq X} \sum_{ 1 \leq p \leq X} \sum_{ b \in \Z} \sum_{ d \in \Z}}_{\substack{ ad -pb =r \\ (a,p) =1}}\alpha(a)  w(b) w(d)\\ 
& = \mathop{\sum_{ 1 \leq a \leq X} \sum_{ 1 \leq p \leq X} \sum_{ b \in \Z }}_{\substack{ b \equiv -r\bar{p} \moda \\ (a,p) = 1}} \alpha(a) w(b) w\l( \frac{r+pb}{a}\r) ,
\ea
where we eliminate $d$ and interpret the equality as a congruence modulo $a$.\\

\begin{proposition} \label{Srzero}
We have,
\be \label{sr0x}
S_r^{(0)}(X)=S_w(X, r)+O\l(HX^{1+\ve}\r).
\ee 
\end{proposition}
\begin{proof}
We first make the easy observation that the total contribution of the terms for which at least one of the variables $b, \frac{r+pb}{a}$ is zero is $O(\tau(r)X)$. \\
We will now look into the contribution of the terms for which at least one of the variables $b,\frac{r+pb}{a}$ is in $(1-H,1)\cup (X,X+H)$. Let us first consider the case when $\frac{r+pb}{a} \in (1-H,1)\cup (X,X+H) $. Then we have,
\ba\label{13} 
\mathop{\sum \sum\sum}_{ \substack{ 1\leq a,p,b \leq X \\ \frac{r+pb}{a} \in  (1-H,1)\cup (X,X+H) \\ pb \equiv - r \moda } } \alpha(a) 
&\ll  X^{\ve}\sum_{1 \leq a \leq X}\sum_{\substack{ \frac{n+r}{a} \in (1-H,1)\cup (X,X+H) \\ n \equiv - r \moda}} \tau(n)\\
& \ll X^{\ve}\sum_{1 \leq a \leq X}\sum_{\substack{ k \in (1-H,1)\cup (X,X+H) }} 1 \\
& \ll HX^{1+\ve}.
\ea
The case when $b$ is in  $(1-H,1)\cup (X,X+H)$ is similar. Therefore, we have
\be
S_r^{(0)}(X)=S_w(X, r)+O\l(HX^{1+\ve}\r).
\ee
\end{proof}
Our task now is to evaluate $ S_w(X,r)$ asymptotically. 
\subsection{Application of Poisson summation}
We  apply the Poisson summation formula (i.e., Lemma \ref{10}) to the sub-sum over  the variable $b$ and
obtain  the following expression from \eqref{defSrwX}.
\ba\label{SW-decomp}
S_w(X, r) 
&= \mathop{\sum_{1 \leq a \leq X} \sum_{1 \leq p \leq X}}_{(a,p) = 1} \alpha(a) \sum_{\substack{b \in \Z \\ b\equiv -r\bar{p}\moda}}w(b)w\l(\frac{r+pb}{a}\r) \\
& =A_w(X, r)+B_w(X, r),
\ea
where we separate the zero frequency and define 
\be\label{Awxr}
A_w(X, r):= \mathop{\sum_{1 \leq a \leq X}\sum_{1 \leq p \leq X}}_{\substack{ (a,p) = 1 }} \frac {\alpha(a)} {a} \widehat{F_{a,p}}\left(0\right),
\ee
and the non-zero frequency is defined by 
\be \label{India1}
B_w(X, r) := \mathop{\sum_{1 \leq a \leq X}\sum_{1 \leq p \leq X}}_{\substack{ (a,p) = 1 }} \frac {\alpha(a)} {a}  \sum_{\substack{n\in \Z\\ n\neq 0}}e\left(\frac{-nr\Bar{p}}{a}\right)\widehat{F_{a,p}}\left(\frac{n}{a}\right),
\ee
where we define for any positive integer $a$ and any integer $c$
\be \label{F(x)}
F_{a,c}(x):=w(x)w\l(\frac{r+cx}{a}\r).
\ee
\subsection{The main term}
In the following proposition we will see that we can extract a main term of size $X^2$ from $A_w(X, r)$. This will be proved in the next section.
\begin{proposition} \label{MT}
We have,
\be \label{awx}
A_w(X,r)= \mathop{\sum_{1 \leq a \leq X}\sum_{1 \leq p \leq X}} \frac {\alpha(a)} {a} \int\limits _{1}^X w\l(\frac{r+px}{a} \r) \diff x + O\l(HX^{1+\ve}\r),
\ee
where the weight function $w$ satisfies the conditions \eqref{w0}-\eqref{w-prime}.
\end{proposition}
\subsection{The error term} The following is the main proposition of this paper where we have used the estimation for the exponential sum over $\bar{p}$ for large moduli, i.e., Lemma \ref{exponential}.
\begin{proposition}\label{Bw}
We have
\ba
B_{w}(X, r)  \ll X^{\frac{7}{4}+\ve}+ \frac{X^{2 + \frac{1}{2}+\ve}}{H} + r^{\frac{1}{2}}\l(\frac{X}{H} \r)^{\frac{3}{2}}X^{\frac{7}{8}+\ve}. 
\ea
\end{proposition}
\subsection{Proof of the theorem}  The main term coming from $S_{-r}^-(X)$ and $S^+_{|r|}(X)$ are the following: 
\ba \label{other1}
\mathop{\sum_{1 \leq a \leq X}\sum_{1 \leq p \leq X}} \frac {\alpha(a)} {a} \int\limits _{1}^X w\l(\frac{-r+px}{a} \r) \diff x + O\l(HX^{1+\ve}\r) ,
\ea 
and 
\ba \label{other2}
\mathop{\sum_{1 \leq a \leq X}\sum_{1 \leq p \leq X}} \frac {\alpha(a)} {a} \int\limits _{1}^X w\l(\frac{|r| - px}{a} \r) \diff x + O\l(HX^{1+\ve}\r). 
\ea 
Note that The error terms corresponding to the other two sums coincide with $B_{w}(X, r)$. We now analyze the error term $B_{w}(X, r)$ in two cases. We first make an observation that for $r \leq HX^{1/4}$, 
\ba \label{compare}
B_{w}(X, r)  \ll X^{\frac{7}{4}+\ve}+ \frac{X^{2 + \frac{1}{2}+\ve}}{H}.
\ea
In that case, the optimal choice of $H$ is 
\ba\label{optimal1}
H = X^{\frac{3}{4}},
\ea
by comparing \eqref{compare} with $O\l( HX^{1+\ve}\r)$ (see \eqref{sr0x}). However, in the complimentary situation,
\ba 
B_{w}(X, r) \ll r^{\frac{1}{2}}\l(\frac{X}{H} \r)^{\frac{3}{2}}X^{\frac{7}{8}+\ve},
\ea
and we make an optimal choice of $H$ to be
\ba \label{optimal2}
H = r^{\frac{1}{5}}X^{\frac{11}{20}}.
\ea
Hence, the theorem follows from Proposition \ref{MT} and Proposition \ref{Bw}, taking into account \eqref{thm}, \eqref{Srx}, \eqref{sr0x}, \eqref{SW-decomp}, \eqref{other1}, \eqref{other2}, \eqref{optimal1} and \eqref{optimal2}.
 
\section{The main term and the proof of Proposition \ref{MT} }
Using \eqref{Awxr}, we write $A_w(X,r)$ as 
\ba
A_w(X,r) = \mathop{\sum_{1 \leq a \leq X}\sum_{1 \leq p \leq X}} \frac {\alpha(a)} {a} \widehat{F_{a,p}}\left(0\right)  - \mathop{\sum_{1 \leq a \leq X}\sum_{1 \leq p \leq X}}_{\substack{ p|a }} \frac {\alpha(a)} {a} \widehat{F_{a,p}}\left(0\right)
\ea 
Now we trivially bound the second term using the properties of the weight function $w$. Thus, we have
\ba 
A_w(X,r) = \mathop{\sum_{1 \leq a \leq X}\sum_{1 \leq p \leq X}} \frac {\alpha(a)} {a} \widehat{F_{a,p}}\left(0\right) + O(X^{1+\ve}).
\ea
Using the properties of the weight function $w$, it follows from \eqref{Awxr}
\ba 
A_w(X,r) &= \mathop{\sum_{1 \leq a \leq X}\sum_{1 \leq p \leq X}} \frac {\alpha(a)} {a} \int\limits _{1}^X w\l(\frac{r+px}{a} \r) \diff x + O\l(HX^{1+\ve}\r).
\ea
Thus we conclude Proposition \ref{MT}.
\section{The error term and the proof of Proposition \ref{Bw}}
We begin analyzing the sum $B_w(X, r)$. For that we first prove two lemmas. One of them limits the size of the length of the dual sum over $n$ in $B_w(X, r)$ and shows that $\widehat{F_{a,p}}\left(\frac{n}{a}\right)$ is negligibly small if $n\gg X^{1+\ve}/H$. The latter provides the bound for the derivative of $\widehat{F_{a, t}}(y)$ for real $t$. 
\begin{lemma}
We have
\be \label{basic}
\widehat{F_{a, c}}(y)\ll_k H^{1-k} y^{-k}\l(1+\left|\frac{c}{a}\right|^{k-1}\r),
\ee
for every integer $k\geq 1$, where $F_{a, c}(y)$ is defined in \eqref{F(x)}. In particular, by taking $k=1$, we have the bound
\be \label{basic2}
\widehat{F_{a, c}}(y)\ll_k  y^{-1}. 
\ee
\end{lemma}

\begin{proof}

Using \eqref{F(x)} we integrate by parts $k$ times, getting
\begin{align*}
\widehat{F_{a, c}}(y) = \frac{1}{(-2\pi iy)^k}\sum_{j=0}^k \int\limits _{-\infty}^{\infty} w^{(j)}(x){\l(\frac{c}{a}\r)}^{k-j}w^{(k-j)}\l(\frac{r+cx}{a}\r)e(-xy)\diff x .
\end{align*}
Note that $w^{(k)}\l(\frac{r+cx}{a}\r)=0$ unless $\frac{r+cx}{a} \in (1-H, 1)\cup (X, X+H)$, and so $w^{(k)}\l(\frac{cx}{a}\r)$ vanishes outside two intervals
of length $O\l(\l|\frac a c\r| H\r)$ each and by \eqref{w-prime}, we obtain the bound $O({\l|\frac {c}{a}\r|}^{k-1}y^{-k}H^{1-k})$ for the term corresponding to $j=0$. 
Now we consider the other summands. Since $w^{(j)}(x)$ is non-zero only in the intervals $(1-H, 1)$ and $(X, X+H)$ for every $j\geq 1$, 
all these terms together  contribute only $O(y^{-k}H^{1-k})$, again by \eqref{w-prime}.
\end{proof}
The next lemma provides an upper bound for the derivative of $ \widehat{F_{a,t}}\left(\frac{n}{a}\right) $ which is balanced with respect to the variables $n,a$ and $t$.
\begin{lemma} \label{derivative}
For real $t$, we have 
\ba 
\l( \widehat{F_{a,t}}\left(\frac{n}{a}\right)  \r)' \ll \frac{aX}{ntH}.
\ea 
\end{lemma}
\begin{proof}
We have
\ba 
\l( \widehat{F_{a,t}}\left(\frac{n}{a}\right)  \r)' = \int \frac{x}{a} w(x) w'\l(\frac{r+tx}{a} \r) e\l(- \frac{nx}{a} \r) \diff x.
\ea
Using integration by parts 
\ba \label{integralestimate}
\l( \widehat{F_{a,t}}\left(\frac{n}{a}\right)  \r)' &= -\frac{1}{2\pi i} \frac{a}{n}\int \frac{1}{a} w(x) w'\l(\frac{r +tx}{a} \r) e\l(- \frac{nx}{a} \r) \diff x \\
& -\frac{1}{2\pi i} \frac{a}{n} \int \frac{x}{a} w'(x) w'\l(\frac{r+tx}{a} \r) e\l(- \frac{nx}{a} \r) \diff x \\
&-\frac{1}{2\pi i} \frac{a}{n} \int \frac{x}{a}\frac{t}{a} w(x) w''\l(\frac{r+tx}{a} \r) e\l(- \frac{nx}{a} \r) \diff x.
\ea
Now we bound the three terms separately.
By the properties of the weight function $w$ (see \eqref{w-support}, \eqref{w-prime}), we know that 
\ba 
w'\l(\frac{r +tx}{a} \r)\neq 0 \iff \frac{r +tx}{a} \in (1-H,1)\cup (X,X+H),
\ea 
and the length of the integral is $\frac{Ha}{t}$. Hence, we have the bound for the first term.
\ba \label{first}
&-\frac{1}{2\pi i} \frac{a}{n}\int \frac{1}{a} w(x) w'\l(\frac{r +tx}{a} \r) e\l(- \frac{nx}{a} \r) \diff x \\
&\ll \frac{a}{n}.\frac{1}{a}.\frac{Ha}{t}.\frac{1}{H} \\
&\ll \frac{a}{nt}. 
\ea
Next, to bound the second term of \eqref{integralestimate}, we use that fact (see \eqref{w-prime}) that the length of the integral is $H$, since
\ba 
w'(x) \neq 0 \iff x \in (1-H,1)\cup (X,X+H),
\ea 
and 
\ba 
\frac{x}{a} \ll \frac{X}{t}
\ea 
from the support of $w'\l(\frac{r+tx}{a} \r)$ (see \eqref{w-support}). Therefore, we have
\ba \label{second}
&-\frac{1}{2\pi i} \frac{a}{n} \int \frac{x}{a} w'(x) w'\l(\frac{r+tx}{a} \r) e\l(- \frac{nx}{a} \r) \diff x \\
&\ll  \frac{a}{n}. \frac{X}{t}.H.\frac{1}{H^2}  \\
&\ll \frac{aX}{ntH}.
\ea
From the support of $w''\l(\frac{r+tx}{a} \r)$ (see \eqref{w-support}, \eqref{w-prime}), we see that
\ba 
\frac{tx}{a} \ll X ,
\ea
and the length of the integral is $\frac{Ha}{t}$. Hence, we have the bound for the last term of \eqref{integralestimate}.
\ba \label{third}
&-\frac{1}{2\pi i} \frac{a}{n} \int \frac{x}{a}\frac{t}{a} w(x) w''\l(\frac{r+tx}{a} \r) e\l(- \frac{nx}{a} \r) \diff x \\
&\ll \frac{a}{n}. \frac{X}{a} .\frac{Ha}{t}. \frac{1}{H^2} \\
&\ll \frac{aX}{ntH}. 
\ea
Combining \eqref{integralestimate}, \eqref{first}, \eqref{second} and \eqref{third}, we prove the lemma.
\end{proof}
We now focus on $B_w(X, r)$. Interchanging the order of summation,
\be \label{India2}
B_w(X, r) := \sum_{1 \leq a \leq X} \frac {\alpha(a)} {a} \sum_{\substack{0 \neq  n \ll \frac{X^{1+\ve}}{H} }}  \sum_{\substack{ 1 \leq p \leq X \\ (a,p) = 1 }}  e\left(\frac{-nr\Bar{p}}{a}\right)\widehat{F_{a,p}}\left(\frac{n}{a}\right).
\ee
We first divide the terms in $B_{w}(X, r)$ according to $(nr,a)= d$ such that $d \geq 1$.
\ba 
B_{w}(X, r) &= \sum_{d \geq 1}\sum_{ 1 \leq a \leq X} \frac {\alpha(a)} {a} \sum_{\substack{0 \neq  n \ll \frac{X^{1+\ve}}{H}  \\  (a,nr) = d } }  \sum_{\substack{ 1 \leq p \leq X \\ (a,p) = 1 }}  e\left(\frac{-nr\Bar{p}}{a}\right)\widehat{F_{a,p}}\left(\frac{n}{a}\right)  + O\l(X ^{-100}\r) \\
& = \sum_{d \geq 1}\frac{1}{d} \sum_{\substack{0 \neq  n \ll \frac{X^{1+\ve}}{H}  \\  d|nr } } \sum_{\substack{ 1 \leq a_1 \leq \frac{X}{d} \\ (a_1, nr/d) =1 }} \frac {\alpha(da_1)} {a_1}\sum_{\substack{ 1 \leq p \leq X \\ (a_1,p) = 1 }}  e\left(\frac{-(nr/d)\Bar{p}}{a_1}\right)\widehat{F_{da_1,p}}\left(\frac{n}{da_1}\right)  .
\ea
We now divide the range of the variable $a_1$ into dyadic intervals and we write 
\ba 
B_w(X, r)  &= \sum_{d \geq 1}\frac{1}{d} \sum_{\substack{ Q \text{ dyadic} \\ Q \leq \frac{X}{d} }} \sum_{\substack{0 \neq  n \ll \frac{X^{1+\ve}}{H}  \\  d|nr } } \sum_{\substack{ a_1 \sim Q \\ (a_1, nr/d) =1  }} \frac {\alpha(da_1)} {a_1}\sum_{\substack{ 1 \leq p \leq X \\ (a_1,p) = 1 }}  e\left(\frac{-(nr/d)\Bar{p}}{a_1}\right)\widehat{F_{da_1,p}}\left(\frac{n}{da_1}\right) .
\ea
Now our plan is to apply Lemma \ref{exponential}. For that we separate the terms in $B_w(X, r)$ for which $Q^{4/3} \geq X$. Thus we write $B_w(X, r)$ as 
\ba 
B_w(X, r) = B_{1,w}(X, r) + B_{2,w}(X, r) ,
\ea
where 
\ba 
B_{1,w}(X, r) := \sum_{d \geq 1}\frac{1}{d} \sum_{\substack{ Q \text{ dyadic} \\ Q \leq X^{3/4}}} \sum_{\substack{0 \neq  n \ll \frac{X^{1+\ve}}{H}  \\  d|nr } } \sum_{\substack{ a_1 \sim Q \\ (a_1, nr/d) =1  }} \frac {\alpha(da_1)} {a_1}\sum_{\substack{ 1 \leq p \leq X \\ (a_1,p) = 1 }}  e\left(\frac{-(nr/d)\Bar{p}}{a_1}\right)\widehat{F_{da_1,p}}\left(\frac{n}{da_1}\right),
\ea 
and 
\ba 
B_{2,w}(X, r) := \sum_{d \geq 1}\frac{1}{d} \sum_{\substack{ Q \text{ dyadic} \\ X^{3/4} < Q \leq \frac{X}{d} }} \sum_{\substack{0 \neq  n \ll \frac{X^{1+\ve}}{H}  \\  d|nr } } \sum_{\substack{ a_1 \sim Q \\ (a_1, nr/d) =1  }} \frac {\alpha(da_1)} {a_1}\sum_{\substack{ 1 \leq p \leq X \\ (a_1,p) = 1 }}  e\left(\frac{-(nr/d)\Bar{p}}{a_1}\right)\widehat{F_{da_1,p}}\left(\frac{n}{da_1}\right),
\ea 
We first bound the sum $B_{1,w}(X, r)$ trivially.
\begin{proposition}\label{B1w}
We have
\ba
B_{1,w}(X, r) \ll X^{7/4+\ve}. 
\ea
\end{proposition}
\begin{proof}
Using \eqref{basic2}, we have
\ba 
B_{1,w}(X, r) &\ll  \sum_{d \geq 1}  \sum_{\substack{ Q \text{ dyadic} \\ Q \leq X^{3/4}}} \sum_{\substack{0 \neq  n \ll \frac{X^{1+\ve}}{H}  \\  d|nr } } \frac {1} {n} \sum_{\substack{ a_1 \sim Q \\ (a_1, nr/d) =1  }}\sum_{\substack{ 1 \leq p \leq X \\ (a_1,p) = 1 }} 1 \\
&\ll X^{7/4 +\ve},
\ea 
since there are $ O\l( X^{\ve}\r)$ many values of $d$ that contribute to the sum over $d$.
\end{proof}
The following proposition is the main result of this section.
\begin{proposition}\label{B2w}
We have
\ba
B_{2,w}(X, r)  \ll \frac{X^{2 + \frac{1}{2}+\ve}}{H} + r^{\frac{1}{2}}\l(\frac{X}{H} \r)^{\frac{3}{2}}X^{\frac{7}{8}+\ve}. 
\ea
\end{proposition}
\begin{proof}
To bound $B_{2,w}(X, r)$, we first use the summation by parts in the sum over $p$ to separate $\widehat{F_{da_1,p}}\left(\frac{n}{da_1}\right)$. Therefore,
\ba \label{boundB2} 
&B_{2,w}(X, r) = \sum_{d \geq 1}\frac{1}{d} \sum_{\substack{ Q \text{ dyadic} \\ X^{3/4} < Q \leq \frac{X}{d} }} \sum_{\substack{0 \neq  n \ll \frac{X^{1+\ve}}{H}  \\  d|nr } } \sum_{\substack{ a_1 \sim Q \\ (a_1, nr/d) =1  }} \frac {\alpha(da_1)} {a_1}\widehat{F_{da_1,X}}\left(\frac{n}{da_1}\right) \sum_{\substack{ 1 \leq p \leq X \\ (a_1,p) = 1 }}  e\left(\frac{-(nr/d)\Bar{p}}{a_1}\right) \\
& + \sum_{d \geq 1}\frac{1}{d} \sum_{\substack{ Q \text{ dyadic} \\ X^{3/4} < Q \leq \frac{X}{d} }} \sum_{\substack{0 \neq  n \ll \frac{X^{1+\ve}}{H}  \\  d|nr } } \sum_{\substack{ a_1 \sim Q \\ (a_1, nr/d) =1  }} \frac {\alpha(da_1)} {a_1} \int\limits_1^X \l(\widehat{F_{da_1,t}}\left(\frac{n}{da_1}\right)\r)' \sum_{\substack{ 1 \leq p \leq t \\ (a_1,p) = 1 }}  e\left(\frac{-(nr/d)\Bar{p}}{a_1}\right) \diff t.
\ea
Next we estimate the first term. Applying \eqref{basic2} and Lemma \ref{exponential},
\ba \label{firstterm} 
&\sum_{d \geq 1}\frac{1}{d} \sum_{\substack{ Q \text{ dyadic} \\ X^{3/4} < Q \leq \frac{X}{d} }} \sum_{\substack{0 \neq  n \ll \frac{X^{1+\ve}}{H}  \\  d|nr } } \sum_{\substack{ a_1 \sim Q \\ (a_1, nr/d) =1  }}\l| \frac {\alpha(da_1)} {a_1}\r|\l|\widehat{F_{da_1,X}}\left(\frac{n}{da_1}\right)\r| \l|\sum_{\substack{ 1 \leq p \leq X \\ (a_1,p) = 1 }}  e\left(\frac{-(nr/d)\Bar{p}}{a_1}\right)\r| \\
&\ll r^{\frac{1}{2}}\l(\frac{X}{H} \r)^{\frac{1}{2}} X^{\frac{7}{8}+\ve}.
\ea 
Now we focus on the second term of \eqref{boundB2}. 
Now, in order to apply Lemma \ref{exponential}, we first take out the integral over $t$ outside the $a_1$-sum and we decompose the $t$ integral into two parts.
\ba \label{last}
&\sum_{d \geq 1}\frac{1}{d} \sum_{\substack{ Q \text{ dyadic} \\ X^{3/4} < Q \leq \frac{X}{d} }} \sum_{\substack{0 \neq  n \ll \frac{X^{1+\ve}}{H}  \\  d|nr } } \int\limits_1^X \sum_{\substack{ a_1 \sim Q \\ (a_1, nr/d) =1  }} \frac {\alpha(da_1)} {a_1}  \l(\widehat{F_{da_1,t}}\left(\frac{n}{da_1}\right)\r)' \sum_{\substack{ 1 \leq p \leq t \\ (a_1,p) = 1 }}  e\left(\frac{-(nr/d)\Bar{p}}{a_1}\right) \diff t \\
&= \sum_{d \geq 1}\frac{1}{d} \sum_{\substack{ Q \text{ dyadic} \\ X^{3/4} < Q \leq \frac{X}{d} }} \sum_{\substack{0 \neq  n \ll \frac{X^{1+\ve}}{H}  \\  d|nr } } \int\limits_1^{X^{1/2}} \sum_{\substack{ a_1 \sim Q \\ (a_1, nr/d) =1  }} \frac {\alpha(da_1)} {a_1}  \l(\widehat{F_{da_1,t}}\left(\frac{n}{da_1}\right)\r)' \sum_{\substack{ 1 \leq p \leq t \\ (a_1,p) = 1 }}  e\left(\frac{-(nr/d)\Bar{p}}{a_1}\right) \diff t \\
&+ \sum_{d \geq 1}\frac{1}{d} \sum_{\substack{ Q \text{ dyadic} \\ X^{3/4} < Q \leq \frac{X}{d} }} \sum_{\substack{0 \neq  n \ll \frac{X^{1+\ve}}{H}  \\  d|nr } } \int\limits_{X^{1/2}}^X \sum_{\substack{ a_1 \sim Q \\ (a_1, nr/d) =1  }} \frac {\alpha(da_1)} {a_1}  \l(\widehat{F_{da_1,t}}\left(\frac{n}{da_1}\right)\r)' \sum_{\substack{ 1 \leq p \leq t \\ (a_1,p) = 1 }}  e\left(\frac{-(nr/d)\Bar{p}}{a_1}\right) \diff t.
\ea
Trivial estimation along with Lemma \ref{derivative}, the first term is bounded by $O\l(\frac{X^{2 + 1/2+\ve}}{H}\r)$. At last we focus on the second term of \eqref{last} and the main contribution will come from this term. We are ready to apply Lemma \ref{exponential} to get a cancellation in our sum which now satisfies the condition $Q^{4/3} \geq t \geq Q^{1/2}$. Hence, by Lemma \ref{derivative} we have
\ba \label{secondterm}
&\sum_{d \geq 1}\frac{1}{d} \sum_{\substack{ Q \text{ dyadic} \\ X^{3/4} < Q \leq \frac{X}{d} }} \sum_{\substack{0 \neq  n \ll \frac{X^{1+\ve}}{H}  \\  d|nr } } \int\limits_{X^{1/2}}^X \sum_{\substack{ a_1 \sim Q \\ (a_1, nr/d) =1  }} \frac {\alpha(da_1)} {a_1}  \l(\widehat{F_{da_1,t}}\left(\frac{n}{da_1}\right)\r)' \sum_{\substack{ 1 \leq p \leq t \\ (a_1,p) = 1 }}  e\left(\frac{-(nr/d)\Bar{p}}{a_1}\right) \diff t \\
&\ll r^{\frac{1}{2}}\l(\frac{X}{H} \r)^{\frac{3}{2}}X^{\frac{7}{8}+\ve}.
\ea
Combining \eqref{boundB2}, \eqref{firstterm}, \eqref{last} and \eqref{secondterm} we arrive at the proposition.
\end{proof}
\subsection{Proof of Proposition \ref{Bw} }
Proposition \ref{Bw} directly follows by bringing together Proposition \ref{B1w} and Proposition \ref{B2w}.

\section{Proof of Theorem \ref{theorem}}
We first consider the case $r = 0$. By symmetry, we have
\ba
\tilde{\mathpzc{S}}_0(X)  &=  8\mathop{\sum\sum\sum\sum}_{\substack{ 1\leq m_1,m_2,p_1,p_2 \leq X\\ m_1p_1 = m_2p_2 \\ (m_1m_2p_1p_2,q) = 1}} 1 \\
& = 8 \sum_{d |q} \mu(d) \mathop{\sum\sum\sum\sum}_{\substack{ 1\leq m_1,m_2 \leq X/d\\ p_1,p_2 \leq X\\ m_1p_1 = m_2p_2 \\ p_1 \neq p_2 \\ (p_1p_2,q) = 1}} 1  + \frac{\phi(q)}{q} \pi(X) X  + O (X^{1+\ve})\\
& =  \sum_{d |q} \mu(d) \mathop{\sum\sum}_{\substack{p_1,p_2 \leq X\\  p_1 \neq p_2 \\ (p_1p_2,q) = 1}} \mathop{\sum}_{\substack{ 1\leq k \leq \min \{X/dp_1, X/dp_2\}}} 1 +  \frac{\phi(q)}{q} \pi(X) X +  O\l(X^{1+\ve}\r),
\ea 
since we take $m_1/p_2= m_2/p_1 = k$ and the terms corresponding to $p_1=p_2$, which automatically imply $m_1=m_2$, give another main term.  By symmetry we have
\ba 
\tilde{\mathpzc{S}}_0(X)  &= 16 \sum_{d |q} \mu(d) \sum_{p_1 \leq X}\sum_{p_2 <p_1}  \l[\frac{X}{dp_1} \r] +  \frac{\phi(q)}{q} \pi(X) X + \frac{\phi(q)}{q} \pi(X) X  +  O\l(X^{1+\ve}\r) \\
& = 16 \sum_{d |q} \mu(d) \sum_{p_1 \leq X} \pi(p_1)  \l[\frac{X}{dp_1} \r] + \frac{\phi(q)}{q} \pi(X) X + O\l( X^{1+\ve}\r) \\
& = 16 \sum_{d |q} \mu(d) \sum_{p_1 \leq X} \sum_{np_1 \leq X/d} \pi(p_1) + \frac{\phi(q)}{q} \pi(X) X  + O\l( X^{1+\ve}\r) .\\
\ea 
Intechanging the sums over $p_1$, $n$ and by the Prime Number Theorem, we have
\ba 
\tilde{\mathpzc{S}}_0(X)  &= 16 \sum_{d |q} \mu(d) \sum_{n \leq X/2d} \sum_{p_1 \leq X/dn} \pi(p)+ \frac{\phi(q)}{q} \pi(X) X + O\l( X^{1+\ve}\r) \\
& = \frac{4\pi^2}{3}
\left(\prod_{p\mid q}\left(1-\frac{1}{p^2}\right)\right)
\frac{X^2}{\log^2 X}  + \frac{\phi(q)}{q} \pi(X) X 
+ O \left(\frac{X^2}{\log^3 X}\right),
\ea 
For $r \neq 0$, we detect the gcd condition by Mobius function and we have,
\ba 
\tilde{\mathpzc{S}}_r(X)  &=  \mathop{\sum\sum\sum\sum}_{\substack{ -X \leq m_1,m_2,p_1,p_2 \leq X\\ m_1p_1 - m_2p_2 = rq \\ (m_1m_1p_1p_2,q) = 1  }} 1 
  = \sum_{\ell |q} \mu(\ell)  \mathop{\sum\sum\sum\sum}_{\substack{ -X/\ell  \leq m_1,m_2 \leq X/\ell \\ p_1,p_2 \leq X\\ m_1p_1 - m_2p_2 = r/ \ell \\ p_1 \neq p_2 \\ (p_1p_2,q) = 1 }} 1 \\
&  =  \sum_{\ell |q} \mu(\ell)  \mathop{\sum\sum\sum\sum}_{\substack{ -X/\ell \leq m_1,m_2 \leq X/\ell \\ p_1,p_2 \leq X\\ m_1p_1 - m_2p_2 = r/ \ell \\ (p_1, p_2) = 1 }} 1 + O\l(  X^{1+\ve }\r).\\
\ea 
Following the proof of Theorem \ref{mainthm}, we arrive at the case of $r \neq 0$ in Theorem \ref{theorem}.
\section{Proof of Theorem \ref{mixedmoment}}
We start with the following sum:
\ba 
\sum_{\chi (\bmod q)} \left| \sum_{m \leq X}  \chi(m)  \right|^2 \left| \sum_{p \leq X}  \chi(p)  \right|^2 &= \phi(q) \mathop{\sum\sum\sum\sum}_{\substack{ 1\leq m_1,m_2,p_1,p_2 \leq X\\ m_1p_1 \equiv m_2p_2(\bmod q)\\ (m_1m_2p_1p_2,q) = 1}} 1 \\
& =  \phi(q) \l( \Sigma_{\text{diagonal}} + \Sigma_{\text{non-diagonal}} \r),
\ea
where 
\ba 
\Sigma_{\text{diagonal}} &= \mathop{\sum\sum\sum\sum}_{\substack{ 1\leq m_1,m_2,p_1,p_2 \leq X\\ m_1p_1 = m_2p_2 \\ (m_1m_2p_1p_2,q) = 1}} 1 ,
\ea 
and 
\ba 
\Sigma_{\text{non-diagonal}} &= \sum_{r \neq 0 } \mathop{\sum\sum\sum\sum}_{\substack{ 1\leq m_1,m_2,p_1,p_2 \leq X\\ m_1p_1 - m_2p_2 = rq \\ (m_1m_2p_1p_2,q) = 1}} 1 .
\ea 
By Theorem \ref{theorem}, the diagonal terms where $m_1p_1= m_2p_2$ contribute the following 
\ba
\Sigma_{\text{diagonal}} & = \frac{\pi^2}{6}
\left(\prod_{p\mid q}\left(1-\frac{1}{p^2}\right)\right)
\frac{X^2}{\log^2 X}  + \frac{\phi(q)}{q} \pi(X) X 
+ O \left(\frac{X^2}{\log^3 X}\right),
\ea 
We now focus on the non-diagonal terms. Following the main term calculation in Theorem \ref{mainthm} and by Theorem \ref{theorem}, we have 
\ba 
\Sigma_{\text{nondiagonal}}  = \sum_{\ell | q} \mu(\ell) \sum_{ r \neq 0 } \mathop{\sum\sum}_{p_1,p_2\leq X} \frac{1}{p_1p_2}\int w\l(\frac{\ell x}{p_2}\r) w\l(\frac{\ell x+ rq}{p_1} \r) \diff x+ O \l(\frac{X^{3+19/20+\ve}}{q} \r),
\ea 
since $rq \ll X^2$. By Poisson Summation Formula on the sum over $r$, we have 
\ba 
\sum_{ r \neq 0}  w\l(\frac{\ell x+ rq}{p_1} \r) &= \sum_{ r }  w\l(\frac{\ell x+ rq}{p_1} \r) + O(1) \\
&= \int w\l(\frac{\ell x+ yq}{p_1} \r) \diff y + O(1),
\ea 
where we use the integration by parts in the integral and the properties of the weight function $w$. Therefore, we have
\ba 
\Sigma_{\text{nondiagonal}}  &= \sum_{\ell | q} \mu(\ell)  \mathop{\sum\sum}_{p_1,p_2\leq X} \frac{1}{p_1p_2}\iint w\l(\frac{\ell x}{p_2}\r) w\l(\frac{\ell x+ yq}{p_1} \r) \diff x \diff y \\
&+ O \l(X^{2+\ve}\r)+ O \l(\frac{X^{3+19/20+\ve}}{q} \r).
\ea 
By change of variable and the properties of the function $w$, we simplify the main term and obtain
\ba 
\Sigma_{\text{nondiagonal}}  = \frac{\phi(q)}{q^2}X^2\pi(X)^2+ O \l(\frac{X^{3+19/20+\ve}}{q} \r)+ O \l(X^{2+\ve}\r).
\ea 
Thus, we have the theorem.
\bibliographystyle{amsalpha}
\bibliography{reference}
 
\end{document}